\documentclass[11pt,twoside]{amsart}

\usepackage{latexsym}
\usepackage{amssymb}

\newcommand{\eps}{\varepsilon}
\newcommand{\de}{\delta}

\newcommand{\N}{{\mathbb N}}

\theoremstyle{plain}
\newtheorem{thm}{Theorem}[section]

\newtheorem{lemma}[thm]{Lemma}

\theoremstyle{definition}

\begin{document}

\title[Coverings by convex bodies and inscribed balls]%
{Coverings by convex bodies and inscribed balls}

\author{Vladimir Kadets}

\date{\today}

\address{Department of Mechanics and Mathematics, 
Kharkov National University,\linebreak
 pl.~Svobody~4,  61077~Kharkov, Ukraine}
\email{vova1kadets@yahoo.com}

\curraddr{Department of Mathematics, University of Missouri,
Columbia MO 65211}

\thanks{I would like to express my thanks to the Department of 
Mathematics, University of Missouri-Columbia and especially to
Professor Nigel Kalton for hospitality and fruitful working
atmosphere.}


\begin{abstract}
Let $H$ be a Hilbert space. For a closed convex 
body $A$ denote by $r(A)$ the supremum of radiuses of balls, contained
in $A$. We prove, that $\sum_{n=1}^\infty r(A_n) \ge r(A)$ for every covering 
of a convex closed body  $A \subset H$ by a sequence of convex closed 
bodies $A_n$, $n \in \N$. It looks like this fact is 
new even for triangles in a 2-dimensional space.
\end{abstract}

\maketitle

\section{Introduction}

Recall, that by a plank of width $w$ in a Hilbert space one means a set of
the form 
$$
P=\{h \in H : |\langle h - h_0 \, , \, e \rangle| \le \frac{w}{2}\},
$$ 
where $\|e\| = 1$.
According to T.Bang's theorem \cite{bang}, if a sequence $P_n$ of planks 
of widths $w_n$ covers a ball of diameter $w$ , then
$\sum w_n \ge w$. 
(To be more precise this is a particular case of Bang's theorem. The complete
statement includes a convex set $P$ of minimal width $w$ instead of a ball).

K. Ball generalized the Bang's theorem to coverings of a ball in a Banach
space, for planks defined with help of linear functionals instead of
inner product. 

Our work is inspired by Bang's and Ball's theorems. The idea was to find
a "symmetric" generalization of this theorems, where the elements of the
covering and the covered set are of the same nature, and they are measured
"from inside". In this paper we prove such a generalization for sets
in a Hilbert space. An analogous question for general Banach spaces
remains open. In the proof we use ideas from Ball's exposition of 
Bang's theorem.

\section{The main result}

For simplicity allover the paper we consider real Hilbert spaces 
(finite- or infinite-dimensional).
All the results are true for complex spaces too: the only thing
which one must change in the proofs is the equation of a real 
hyperplane: in complex case such an equation uses not the inner 
product itself, but the real part of inner product.

Let $H$ be a Hilbert space, and let $B$ and $S$ be
its unit ball and its unit sphere respectively. 
By "body" in $H$ we mean a closed subset with non-empty
interior. For a convex body $A$ denote by $r(A)$ 
the supremum of radiuses of balls $x+aB$, contained in $A$. 

If $A$ is a ball of radius $r$, then $r(A)=r$. If $A$ is
a plank of width $h$, then $r(A)=\frac{h}{2}$. 

\begin{thm} \label{main}
If a convex body  $A \subset H$  is covered by a sequence 
of convex bodies $A_n$, $n \in \N$, then 
$\sum_{n=1}^\infty r(A_n) \ge r(A)$.
\end{thm}

To prove the theorem we need first some lemmas. 
The goal of the lemmas is to prove, that a convex 
bounded body in a Hilbert space can be approximated
in some sense by a polytope with finite number of faces. 

\begin{lemma} \label{ball1}
Let $A \subset H$ be a bounded convex body. 
Then there is a ball $U \subset A$ with $r(U)=r(A)$.
\end{lemma}
\begin{proof}
Denote $r(A)=r$. For every $p \in [0,r)$ consider the set $A_p$,
consisting of all those $x \in A$, for which $x+pB \subset A$.
Each $A_p$ is a bounded convex closed set, and hence it is a weak compact.
Since $A_p \neq \emptyset$ and decrease as $p \to r$, there is a
point $o_1 \in \bigcap_{p \in [0,r)} A_p$. Then $U=o_1 + rB$ will be the ball
we need.
\end{proof}

\begin{lemma} \label{comb1}
Let $A \subset H$ be a bounded convex body, $r(A)=r$, and let 
$U=o_1 + rB$ be the ball from the previous lemma. Then for every
$\eps > 0$ there is a finite subset 
$G=\{g_1, \ldots , g_n\} \subset rB$ such that 
$(1+\eps) g_i + o_1 \in H \setminus A$ and 
\begin{equation} \label{conv1}
{\rm dist}\,({\rm conv}\{g_i\}_{i=1}^n, 0) \le \eps.
\end{equation}

\end{lemma}
\begin{proof}
Without loss of generality we may assume $o_1=0$ (otherwise
shift the picture). We must prove, that the closed convex hull of 
$(1+\eps)rB \setminus A$ contains 0.

Assume it is not so. By the Hahn-Banach theorem, 
there is a hyperplane $P$, strictly separating $(1+\eps)rB \setminus A$ 
from 0. This means, that the set $F$ - the bigger part 
of the ball $(1+\eps)rB$, 
lying on the same side of $P$ as 0, is included in $A$. 
The subset $F \bigcup rB$ of $A$ evidently contains a ball
of radius bigger than $r$, which contradicts the definition
of $r$.
\end{proof}

\begin{lemma} \label{separ}
Let $x,y \in H$, $\|x\|= 1 + \eps$ and the hyperplane 
$$
\{h \in H: \langle h, y \rangle = \| y \|^2\}
$$
generated by $y$ separates $x$ from the unit ball.
Then $\|x-y\| \le \de(\eps)$, where $\de(\eps)$ tends to 0
as $\eps$ tends to 0.
\end{lemma}
\begin{proof}
Our conditions mean that 
$$
\langle x, y \rangle \ge \|y\|^2
$$
and
$$
\|y\| = {\rm sup}_{h \in B} \langle h, y \rangle \le \|y\|^2,
$$
i.e. $\|y\| \ge 1.$
We have
$$
\|x-y\|^2 = \|x\|^2 - \|y\|^2 + 2(\|y\|^2 - \langle x, y \rangle)
\le \|x\|^2 - \|y\|^2 \le (1+\eps)^2 - 1.
$$
\end{proof}

\begin{lemma} \label{polytop1}
Let $A \subset H$ be a bounded convex body, $r(A)=r$, and let 
$U=o_1 + rB$ be the ball from lemma \ref{ball1}. Then for every
$\eps > 0$ there is a finite subset 
$V=\{v_1, \ldots , v_n\} \subset rB$  ($n$ depends on $\eps)$)
such that
\begin{equation} \label{polytop2}
A \subset W := \bigcap_{v \in V} 
\{ h \in H: \langle h - o_1, v \rangle < r^2 + \eps \}, \  and
\end{equation}
\begin{equation} \label{conv2}
{\rm dist}\,({\rm conv} V, 0) \le \eps.
\end{equation}
\end{lemma}
\begin{proof}
Apply lemma \ref{comb1} for a small $\eps_1$, to obtain corresponding
$G=\{g_1, \ldots , g_n\} \subset rB$. Since 
$(1+\eps_1) g_i + o_1 \in H \setminus A$ , one can separate
$(1+\eps_1) g_i + o_1$ from $A$ by a hyperplane 
$$
P_i = \{h \in H: \langle h - o_1, w_i \rangle = \| w_i \|^2\}.
$$
Since $U=o_1 + rB \subset A$, $P_i$ separates $U$ from  $(1+\eps_1) g_i + o_1$.
This means that hyperplane 
$$
\{h \in H: \langle h, w_i \rangle = \| w_i \|^2\}
$$
separates $rB$ from $(1+\eps_1) g_i$, and by lemma \ref{separ}
$w_i$ is "very close" to $g_i$, and $\|w_i\|$ is "very close" to $r$.
This means in turn, that when $\eps_1$ is small enough, the
elements  $v_k=r \frac{w_k}{\|w_k\|}$, $k=1,2, \ldots , n$ 
fulfill the demands of this lemma.
\end{proof}

\begin{lemma} \label{polytop}
For a convex body $A \subset B$ and for arbitrary $\de > 0$ there is a
polytope $W$ of the form 
$$\bigcap_{k=1}^n 
\{ h \in H: \langle h - o_1, v_k \rangle < a_k\},
$$
such that $W \supset A$ and $r(W \bigcap B) \le r(A) + \de$
\end{lemma}
\begin{proof}
The polytope $W$ can be taken from (\ref{polytop2}) of the previous 
lemma \ref{polytop1}, with $\eps$ small enough. Let us show this.
Consider $U=o_1 + rB$ from lemma \ref{polytop1}. For simplicity
assume $o_1 = 0$ (the general case differs not too much from this one).
Assume, contrary to our statement, that $W \bigcap B$ contains a ball
of the form $U_0=o_0 + (r+\de)B$ (and automatically $\|o_0\| \le 1$).
Then $W \supset{\rm conv}\{U_0,U\}$. According to condition 
(\ref{conv2}), since $(1+\eps)V$ lies outside $W$,
$$
{\rm dist}\,\left({\rm conv} ((1+\eps)U \setminus W) \,,\, 0  \right) 
\le (1+\eps)\eps,
$$
and hence
$$
{\rm dist}\,\left ({\rm conv} ((1+\eps)U \setminus {\rm conv}\{U_0,U\}) \,,\, 0 \right ) 
\le (1+\eps)\eps.
$$
But for fixed $\de$ and $\eps \to 0$ the last inequality cannot be true,
since radius of $U_0$ equals $(r+\de)$ and the distance between the centers
of $U_0$ and $U$ is bounded by a number, independent on $\eps$.
\end{proof}

\vspace{5 mm}
{\bf Proof of the theorem \ref{main}.}
Denote $r(A_n)=r_n$. By the definition of $r(A)$ we must prove that
$\sum_{n=1}^N r_n \ge r(U)$ for every ball $U \subset A$.  
By homogeneity this means that it is sufficient to prove the theorem
for the case of $A$ being the unit ball $B$ of $H$. Also without loss of 
generality one may assume $A_n \subset B$: otherwise consider 
$A_n \bigcap B$ instead of $A_n$. 

So the theorem can be reformulated as
follows: 

\vspace{5 mm}

Let $\sum_{n=1}^\infty r_n < 1$ and let $A_n \subset B$ be convex bodies
with $r(A_n)=r_n$. Then 
$$
B \setminus \bigcup_{n=1}^\infty A_n \neq \emptyset.
$$ 

\vspace{5 mm}

According to lemma \ref{polytop} every $A_n$ may be included
in interior of a polytope (say, $C_n$), in such  a way, that the
condition $\sum_{n=1}^\infty r(C_n \bigcap B) < 1$ still takes place.
So, in fact, we may assume that interiors of all of $A_n$ are
relatively weakly open in $B$, and it is enough to prove, that $B$ 
cannot be covered by interiors of $A_n$. But in this case, since
$B$ is a weak compact, it is sufficient to prove, that $B$ 
cannot be covered by union of finitely many interiors of $A_n$.
So we reduced our theorem to the case, when the number of sets
is finite (say, equals $N$). 

Let us fix a sequence of positive $\de_n$ 
with 
\begin{equation} \label{delta}
\sum_{n=1}^N (1+ \de_n)r_n < 1,
\end{equation}
and select $\eps_n > 0$ with 
\begin{equation} \label{eps}
(1+ \de_n)(r_n^2 - 6 N \eps_n) \ge r_n^2 + \eps_n
\end{equation}

For every $A_n$ apply lemma \ref{polytop1} with $\eps =\eps_n$
to get corresponding vectors $o_n \in A_n$, finite sets 
$V_n \subset r_n B$ and open polytops 
$$
W_n := \bigcap_{v \in V} 
\{ h \in H: \langle h - o_n, v \rangle < r_n^2 + \eps_n \}, 
$$
such that $A_n \subset W_n$ and
\begin{equation} \label{conv6}
{\rm dist}\,({\rm conv} V_n, 0) \le \eps_n.
\end{equation}

We are going to prove that 

\begin{equation} \label{wow}
B \setminus \bigcup_{n=1}^N W_n \neq \emptyset.
\end{equation}

Let us introduce an auxiliary space 
$X=H \oplus H_1 \oplus H_2 \oplus \ldots \oplus H_N$ -
the Hilbertian orthogonal direct sum of isometric copies of the
original space $H$. Let $U_n: H \to H_n$ be corresponding isometries.

Consider $K=V_1 \times V_2 \times \ldots \times V_N$ and for every 
$\bar{g}=(g_1,g_2, \ldots , g_N) \in K$ introduce
$$
f(\bar{g})= \sum_{n=1}^N ((1 + \de_n)(g_n + U_n g_n) - U_n o_n). 
$$
Fix $\bar{x}=(x_1,x_2, \ldots , x_N) \in K$ for which $\|f(\bar{x})\|$
is maximal.  We are going to prove that 
$$
x= \sum_{n=1}^N (1 + \de_n)x_n
$$
is the element, we need: namely, 
$x$ belongs $B$ but does not belong to any
of $W_j$.

First of all, due to (\ref{delta})
$\| \sum_{n=1}^N (1 + \de_n)x_n\| \le \sum_{n=1}^N (1 + \de_n)r_n < 1$,
i.e. $x \in B$. Now for a fixed $j \in \{1,2, \ldots , N \}$ consider
$$
y_j=f(\bar{x}) - (1 + \de_j)(x_j + U_j x_j).
$$ 
Taking in account,
that according to our construction 
$$
\|y_j + (1 + \de_j)(x_j + U_j x_j)\|=\|f(\bar{x})\| \ge 
\|y_j + (1 + \de_j)(v + U_j v)\|
$$
for all $v \in V_j$, we obtain for all  $v \in V_j$
\begin{equation} \label{sic}
\langle y_j \, , \, (1 + \de_j)(x_j + U_j x_j) \rangle \ge 
\langle y_j \, , \, (1 + \de_j)(v + U_j v) \rangle.
\end{equation} 

Due to condition (\ref{conv6}), there is a convex combination
$\sum_{v \in V_j} \alpha_v v$, having norm less then $\eps_j$.
Multiplying (\ref{sic}) by $\alpha_v$ and adding for all $v \in V_j$
we obtain
$$
\langle y_j \, , \, (1 + \de_j)(x_j + U_j x_j) \rangle \ge
-2 \eps_j (1 + \de_j) \|y_j\| \ge -6 N \eps_j (1 + \de_j)  
$$

Let us transform the last inequality using definition of
$y_j$ and pairwise orthogonality of subspaces $H, H_1, \ldots , H_N$.

$$
\left \langle x - (1 + \de_j)x_j + (1 + \de_j)\sum_{n \neq j} U_n x_n - 
\sum_{n=1}^N U_n o_n \, \,
, \, \, 
 (1 + \de_j)(x_j + U_j x_j ) \right \rangle \ge -6 N \eps_j (1 + \de_j),
$$ 

$$
(1 + \de_j)\langle x \, , \, x_j \rangle - (1 + \de_j)^2 r_j^2
- (1 + \de_j)\langle U_j o_j \, , \, U_j x_j \rangle \ge 
-6 N \eps_j (1 + \de_j),
$$

$$
(1 + \de_j)\langle x \, , \, x_j \rangle  
- (1 + \de_j)\langle o_j \, , \, x_j \rangle \ge 
(1 + \de_j)^2 r_j^2 - 6 N \eps_j (1 + \de_j).
$$
By (\ref{eps}) we have 
$$
\langle x - o_j\, , \, x_j \rangle \ge r_j^2 + \eps_j,
$$
and this by definition of $W_j$ means that $x$ does not belong to $W_j$.

\mbox

\end{document}